\documentclass[12]{article}
\usepackage{graphicx}
\usepackage[top=1in, bottom=1in, left=1in, right=1in]{geometry}
\newtheorem{theorem}{Theorem}
\usepackage{color}
\usepackage{amsmath}

\begin{document}
\title{A Non-Sieving Application of the Euler Zeta Function}
\author{Michael P. May}
\date{October 2, 2015}
\maketitle
\begin{center}
\textit{To Nomi Bell May, my beloved wife and eternal companion...}
\end{center}

\noindent The reader is invited to recall the celebrated equation \cite{HowEulerDiscoveredthezetafunction}

\begin{equation}
\zeta(s) = 1 +\frac{1}{2^s}+\frac{1}{3^s}+\frac{1}{4^s}+\frac{1}{5^s}+ \cdots = \frac{1}{1-\frac{1}{2^s}}\cdot\frac{1}{1-\frac{1}{3^s}}\cdot\frac{1}{1-\frac{1}{5^s}}\cdot\frac{1}{1-\frac{1}{7^s}}\cdot\cdots
\end{equation}

\noindent known as the Euler zeta function.  One who is familiar with this famous equation, discovered by Leonhard Euler, which established a special relationship between the prime and composite numbers might naturally ponder the results of the application of this special function in cases where there is no known way to sieve composite numbers out of the product term in the denominator of the right-hand side of the equation.  Such might be the case when an infinite series of numbers to be analyzed by this function are calculated by a polynomial expression that yields successively increasing positive integer values, though this non-sieving approach could also be applied to the realm of positive integers that Euler used to sieve out the composite terms to arrive at this monumental result.  But when one replaces the positive integers in the original Euler zeta function with successive numbers calculated by an integer-valued polynomial function which will yield prime numbers when the input domain is the positive integers themselves, there may not be an intuitive way to eliminate the composite terms from the denominator on the right-hand side of Eq. 1 by either scaling a previous prime number calculation or by employing predictable values of the domain of the function which would make future outputs of the polynomial prime.  So the best one may be able to hope for in this case is to determine some value to be added or subtracted from unity in the numerator on the right-hand side of Eq. 1 to make both sides equal in the hope that we can analytically use that value to predict the number of prime numbers that exist as outputs of the polynomial function for some limit $x$ of the input domain.  To this end, the following equation is introduced as an application of a non-sieving approach to solving the Euler zeta function for the case of integer-valued polynomials in general:

\newcommand{\opM}{\mathop{\vphantom{\sum}\mathchoice
  {\vcenter{\hbox{\huge M}}}
  {\vcenter{\hbox{\Large A}}}{\mathrm{A}}{\mathrm{A}}}\displaylimits}

\begin{equation}
\scalebox{1.5}{Z(s)} = \displaystyle 1 + \displaystyle\sum\limits_{n=1}^\infty \frac{1}{f(n)^s} = \frac{1+\displaystyle\opM\limits_{i=1}^\infty \Big[(-1)^{i}\Big]\displaystyle\sum{}^{i+1}}{\displaystyle\prod\limits_{n=1}^\infty {1-\frac{1}{\displaystyle f(n)^s}•}}
\end{equation}

\noindent where $f(n)$ is some integer-valued polynomial function that has as its domain the realm of positive integers and $f(1) > 1$.  If $f(1) = 1$, then Eq. 2 becomes

\begin{equation}
\scalebox{1.5}{Z(s)} = \displaystyle \sum\limits_{n=1}^\infty \frac{1}{f(n)^s} = \frac{1+\displaystyle\opM\limits_{i=1}^\infty \Big[(-1)^{i}\Big]\displaystyle\sum{}^{i+1}}{\displaystyle\prod\limits_{n=2}^\infty {1-\frac{1}{\displaystyle f(n)^s}•}}
\end{equation}

A definition of the "M-series" term in these big Zeta equations is in order to lay the foundation for the understanding of the application of Euler's zeta function to the output of integer-valued polynomials using the non-sieving approach:

$$\opM\limits_{i=2}^\infty \Big[(-1)^{i-1}\Big]\displaystyle\sum{}^i=-\Sigma_1 \Sigma_2+\Sigma_1 \Sigma_2 \Sigma_3-\Sigma_1 \Sigma_2 \Sigma_3 \Sigma_4+\Sigma_1 \Sigma_2 \Sigma_3 \Sigma_4 \Sigma_5-\Sigma_1 \Sigma_2 \Sigma_3 \Sigma_4 \Sigma_5 \Sigma_6+...$$

\noindent where\\

$$\Sigma_1 \Sigma_2 = \sum\limits_{i=2}^\infty \sum\limits_{j=i}^\infty \frac{1}{f(i)} \frac{1}{f(j)}$$

$$\Sigma_1 \Sigma_2 \Sigma_3 = \sum\limits_{i=2}^\infty \sum\limits_{j=i}^\infty \sum\limits_{k=j+1}^\infty \frac{1}{f(i)} \frac{1}{f(j)} \frac{1}{f(k)}$$

$$\Sigma_1 \Sigma_2 \Sigma_3 \Sigma_4 = \sum\limits_{i=2}^\infty \sum\limits_{j=i}^\infty \sum\limits_{k=j+1}^\infty \sum\limits_{l=k+1}^\infty \frac{1}{f(i)} \frac{1}{f(j)} \frac{1}{f(k)} \frac{1}{f(l)}$$

$$\Sigma_1 \Sigma_2 \Sigma_3 \Sigma_4 \Sigma_5 = \sum\limits_{i=2}^\infty \sum\limits_{j=i}^\infty \sum\limits_{k=j+1}^\infty \sum\limits_{l=k+1}^\infty \sum\limits_{m=l+1}^\infty \frac{1}{f(i)} \frac{1}{f(j)} \frac{1}{f(k)} \frac{1}{f(l)}\frac{1}{f(m)}$$

$$\Sigma_1 \Sigma_2 \Sigma_3 \Sigma_4 \Sigma_5 \Sigma_6= \sum\limits_{i=2}^\infty \sum\limits_{j=i}^\infty \sum\limits_{k=j+1}^\infty \sum\limits_{l=k+1}^\infty \sum\limits_{m=l+1}^\infty \sum\limits_{n=m+1}^\infty \frac{1}{f(i)} \frac{1}{f(j)} \frac{1}{f(k)} \frac{1}{f(l)} \frac{1}{f(m)} \frac{1}{f(n)}$$

\begin{center}
\bf.
\end{center}
\begin{center}
\bf.
\end{center}
\begin{center}
\bf.
\end{center}

\noindent As an example to support the notion of a non-sieving application of the Euler zeta function to integer-valued polynomial functions, let's apply this non-sieving approach to the polynomial $$n^3-(n-1)^3 = 3n^2-3n+1$$ for which the domain of the base $n$ is the positive integers less than or equal to some limit $x$ that generate successive integer values, prime and composite, for processing in the Z function.  For convenience in notation and to generalize the notion, the polynomial function $n^2-3n+1$ will be represented by $f(n)$.  This polynomial is defined in another one of the author's research papers as a prime shell of degree three, one of a family of integer-valued functions which are analyzed using this non-sieving approach\cite{On the Existence and Frequency of the Shell Primes}.\\  

\noindent We begin our example of the non-sieving application of the Euler zeta function to the polynomial $n^2-3n+1$ by writing

$$Z(s) = \frac{1}{f(1)^s} +\frac{1}{f(2)^s}+\frac{1}{f(3)^s}+\frac{1}{f(4)^s}+\frac{1}{f(5)^s}+\frac{1}{f(6)^s}+\frac{1}{f(7)^s}+\frac{1}{f(8)^s}+\frac{1}{f(9)^s}+\frac{1}{f(10)^s}+...$$

\noindent or,

$$Z(s) = 1 +\frac{1}{7^s}+\frac{1}{19^s}+\frac{1}{37^s}+\frac{1}{61^s}+\frac{1}{91^s}+\frac{1}{127^s}+\frac{1}{169^s}+\frac{1}{217^s}+\frac{1}{271^s}+...$$

\noindent for which we will drop the powers of $s$ for convenience in notation, i.e.,

\begin{equation}
Z = 1 +\frac{1}{7}+\frac{1}{19}+\frac{1}{37}+\frac{1}{61}+\frac{1}{91}+\frac{1}{127}+\frac{1}{169}+\frac{1}{217}+\frac{1}{271}+....
\end{equation}\\

\noindent We first multiply both sides of Eq. 4 by the first fraction on the right hand side of that equation which yields

$$\frac{1}{7}\cdot Z =\frac{1}{7}+\frac{1}{7}\cdot\frac{1}{7}+\frac{1}{7}\cdot\frac{1}{19}+\frac{1}{7}\cdot\frac{1}{37}+\frac{1}{7}\cdot\frac{1}{61}+\frac{1}{7}\cdot\frac{1}{91}+\frac{1}{7}\cdot\frac{1}{127}+\frac{1}{7}\cdot\frac{1}{169}+\frac{1}{7}\cdot\frac{1}{217}+\frac{1}{7}\cdot\frac{1}{271}+...$$

\noindent and when we subtract this result from Eq. 4, we get

\begin{equation}
\Big(1-\frac{1}{7}\Big)\cdot Z =1 +\frac{1}{19}+\frac{1}{37}+\frac{1}{61}+\frac{1}{91}+\frac{1}{127}+\frac{1}{169}+\frac{1}{217}+\frac{1}{271}+...-\displaystyle\sum\limits_{n=2}^\infty \frac{1}{7 \cdot f(n)}.
\end{equation}

\noindent Continuing with the multiplication of this new equation by the next fraction in the sequence, $\frac{1}{19}$, 

\begin{multline*}
$$\frac{1}{19}\cdot\Big(1-\frac{1}{7}\Big)\cdot Z =\frac{1}{19}+\frac{1}{19}\cdot\frac{1}{19}+\frac{1}{19}\cdot\frac{1}{37}+\frac{1}{19}\cdot\frac{1}{61}+\frac{1}{19}\cdot\frac{1}{91}+\frac{1}{19}\cdot\frac{1}{127}+\frac{1}{19}\cdot\frac{1}{169}+\frac{1}{19}\cdot\frac{1}{217}+\frac{1}{19}\cdot\frac{1}{271}+...\\-\displaystyle\sum\limits_{n=2}^\infty \frac{1}{7 \cdot 19\cdot f(n)}$$
\end{multline*}

\noindent so that when we subtract this result from Eq. 5, we end up with

\begin{multline*}
$$\Big(1-\frac{1}{19}\Big)\Big(1-\frac{1}{7}\Big)\cdot Z =1 +\frac{1}{37}+\frac{1}{61}+\frac{1}{91}+\frac{1}{127}+\frac{1}{169}+\frac{1}{217}+\frac{1}{271}+...\\-\displaystyle\sum\limits_{n=2}^\infty \frac{1}{7 \cdot f(n)}-\displaystyle\sum\limits_{n=3}^\infty \frac{1}{19 \cdot f(n)}+\displaystyle\sum\limits_{n=2}^\infty \frac{1}{7\cdot19 \cdot f(n)}.$$
\end{multline*}

\noindent Continuing with the multiplication of this new equation by the next term in the sequence, $\frac{1}{37}$,

\begin{multline*}
$$\frac{1}{37}\cdot\Big(1-\frac{1}{19}\Big)\cdot\Big(1-\frac{1}{7}\Big)\cdot Z =\frac{1}{37}+\frac{1}{37}\cdot\frac{1}{37}+\frac{1}{37}\cdot\frac{1}{61}+\frac{1}{37}\cdot\frac{1}{91}+\frac{1}{37}\cdot\frac{1}{127}+\frac{1}{37}\cdot\frac{1}{169}+\frac{1}{37}\cdot\frac{1}{217}+\frac{1}{37}\cdot\frac{1}{271}+...\\-\displaystyle\sum\limits_{n=2}^\infty \frac{1}{7 \cdot37\cdot f(n)}-\displaystyle\sum\limits_{n=3}^\infty \frac{1}{19 \cdot37\cdot f(n)}+\displaystyle\sum\limits_{n=2}^\infty \frac{1}{7 \cdot 19\cdot 37\cdot f(n)}$$
\end{multline*}

\noindent so that

\begin{multline*}
$$\Big(1-\frac{1}{37}\Big)\Big(1-\frac{1}{19}\Big)\Big(1-\frac{1}{7}\Big)\cdot Z =1 +\frac{1}{61}+\frac{1}{91}+\frac{1}{127}+\frac{1}{169}+\frac{1}{217}+\frac{1}{271}+...\\-\displaystyle\sum\limits_{n=2}^\infty \frac{1}{7 \cdot f(n)}-\displaystyle\sum\limits_{n=3}^\infty \frac{1}{19 \cdot f(n)}-\displaystyle\sum\limits_{n=4}^\infty \frac{1}{37 \cdot f(n)}\\+\displaystyle\sum\limits_{n=2}^\infty \frac{1}{7\cdot19 \cdot f(n)}+\displaystyle\sum\limits_{n=2}^\infty \frac{1}{7 \cdot37\cdot f(n)}+\displaystyle\sum\limits_{n=3}^\infty \frac{1}{19 \cdot37\cdot f(n)}-\displaystyle\sum\limits_{n=2}^\infty \frac{1}{7 \cdot 19\cdot 37\cdot f(n)}.$$
\end{multline*}

\noindent Continuing with the multiplication of this new equation by the next term in the sequence, $\frac{1}{61}$,

\begin{multline*}
$$\frac{1}{61}\cdot\Big(1-\frac{1}{37}\Big)\cdot\Big(1-\frac{1}{19}\Big)\cdot\Big(1-\frac{1}{7}\Big)\cdot Z =\frac{1}{61}+\frac{1}{61}\cdot\frac{1}{61}+\frac{1}{61}\cdot\frac{1}{91}+\frac{1}{61}\cdot\frac{1}{127}+\frac{1}{61}\cdot\frac{1}{169}+\frac{1}{61}\cdot\frac{1}{217}+\frac{1}{61}\cdot\frac{1}{271}+...\\-\displaystyle\sum\limits_{n=2}^\infty \frac{1}{7\cdot 61 \cdot f(n)}-\displaystyle\sum\limits_{n=3}^\infty \frac{1}{19 \cdot61\cdot f(n)}-\displaystyle\sum\limits_{n=4}^\infty \frac{1}{37 \cdot 61\cdot f(n)}+\displaystyle\sum\limits_{n=2}^\infty \frac{1}{7\cdot19 \cdot61\cdot f(n)}\\+\displaystyle\sum\limits_{n=2}^\infty \frac{1}{7 \cdot37\cdot61\cdot f(n)}+\displaystyle\sum\limits_{n=3}^\infty \frac{1}{19 \cdot37\cdot61\cdot f(n)}-\displaystyle\sum\limits_{n=2}^\infty \frac{1}{7 \cdot 19\cdot 37\cdot61\cdot f(n)}$$
\end{multline*}

\noindent so that

\begin{multline*}
$$\Big(1-\frac{1}{61}\Big)\Big(1-\frac{1}{37}\Big)\Big(1-\frac{1}{19}\Big)\Big(1-\frac{1}{7}\Big)\cdot Z =1 +\frac{1}{91}+\frac{1}{127}+\frac{1}{169}+\frac{1}{217}+\frac{1}{271}+...\\-\displaystyle\sum\limits_{n=2}^\infty \frac{1}{7 \cdot f(n)}-\displaystyle\sum\limits_{n=3}^\infty \frac{1}{19 \cdot f(n)}-\displaystyle\sum\limits_{n=4}^\infty \frac{1}{37 \cdot f(n)}-\displaystyle\sum\limits_{n=5}^\infty \frac{1}{61 \cdot f(n)}\\+\displaystyle\sum\limits_{n=2}^\infty \frac{1}{7\cdot19 \cdot f(n)}+\displaystyle\sum\limits_{n=2}^\infty \frac{1}{7 \cdot37\cdot f(n)}+\displaystyle\sum\limits_{n=3}^\infty \frac{1}{19 \cdot37\cdot f(n)}-\displaystyle\sum\limits_{n=2}^\infty \frac{1}{7 \cdot 19\cdot 37\cdot f(n)}+\displaystyle\sum\limits_{n=2}^\infty \frac{1}{7\cdot 61 \cdot f(n)}+\displaystyle\sum\limits_{n=3}^\infty \frac{1}{19 \cdot61\cdot f(n)}\\+\displaystyle\sum\limits_{n=4}^\infty \frac{1}{37 \cdot 61\cdot f(n)}-\displaystyle\sum\limits_{n=2}^\infty \frac{1}{7\cdot19 \cdot61\cdot f(n)}-\displaystyle\sum\limits_{n=2}^\infty \frac{1}{7 \cdot37\cdot61\cdot f(n)}-\displaystyle\sum\limits_{n=3}^\infty \frac{1}{19 \cdot37\cdot61\cdot f(n)}+\displaystyle\sum\limits_{n=2}^\infty \frac{1}{7 \cdot 19\cdot 37\cdot61\cdot f(n)}.$$
\end{multline*}

\noindent Continuing with the multiplication of of this new equation by the next term in the sequence, $\frac{1}{91}$,

\begin{multline*}
$$\frac{1}{91}\cdot\Big(1-\frac{1}{61}\Big)\cdot\Big(1-\frac{1}{37}\Big)\cdot\Big(1-\frac{1}{19}\Big)\cdot\Big(1-\frac{1}{7}\Big)\cdot Z =\frac{1}{91}+\frac{1}{91}\cdot\frac{1}{91}+\frac{1}{91}\cdot\frac{1}{127}+\frac{1}{91}\cdot\frac{1}{169}+\frac{1}{91}\cdot \frac{1}{217}+\frac{1}{91}\cdot \frac{1}{271}+...\\-\displaystyle\sum\limits_{n=2}^\infty \frac{1}{7 \cdot 91 \cdot f(n)}-\displaystyle\sum\limits_{n=3}^\infty \frac{1}{19 \cdot91 \cdot f(n)}-\displaystyle\sum\limits_{n=4}^\infty \frac{1}{37 \cdot 91\cdot f(n)}-\displaystyle\sum\limits_{n=5}^\infty \frac{1}{61 \cdot91 \cdot f(n)}+\displaystyle\sum\limits_{n=2}^\infty \frac{1}{7\cdot19 \cdot 91\cdot f(n)}\\+\displaystyle\sum\limits_{n=2}^\infty \frac{1}{7 \cdot37\cdot91\cdot f(n)}+\displaystyle\sum\limits_{n=3}^\infty \frac{1}{19 \cdot37\cdot91 \cdot f(n)}-\displaystyle\sum\limits_{n=2}^\infty \frac{1}{7 \cdot 19\cdot 37 \cdot 91 \cdot f(n)}+\displaystyle\sum\limits_{n=2}^\infty \frac{1}{7 \cdot 61 \cdot 91 \cdot f(n)}+\displaystyle\sum\limits_{n=3}^\infty \frac{1}{19 \cdot 61 \cdot91 \cdot f(n)}\\+\displaystyle\sum\limits_{n=4}^\infty \frac{1}{37 \cdot 61 \cdot 91 \cdot f(n)}-\displaystyle\sum\limits_{n=2}^\infty \frac{1}{7 \cdot 19 \cdot 61\cdot 91 \cdot f(n)}-\displaystyle\sum\limits_{n=2}^\infty \frac{1}{7 \cdot 37 \cdot 61 \cdot 91 \cdot f(n)}-\displaystyle\sum\limits_{n=3}^\infty \frac{1}{19 \cdot 37 \cdot 61 \cdot 91 \cdot f(n)}\\+\displaystyle\sum\limits_{n=2}^\infty \frac{1}{7 \cdot 19 \cdot 37 \cdot 61 \cdot 91 \cdot f(n)}$$
\end{multline*}

\noindent so that

\begin{multline*}
$$\Big(1-\frac{1}{91}\Big)\Big(1-\frac{1}{61}\Big)\Big(1-\frac{1}{37}\Big)\Big(1-\frac{1}{19}\Big)\Big(1-\frac{1}{7}\Big)\cdot Z =1 +\frac{1}{127}+\frac{1}{169}+\frac{1}{217}+\frac{1}{271}+...\\-\displaystyle\sum\limits_{n=2}^\infty \frac{1}{7 \cdot f(n)}-\displaystyle\sum\limits_{n=3}^\infty \frac{1}{19 \cdot f(n)}-\displaystyle\sum\limits_{n=4}^\infty \frac{1}{37 \cdot f(n)}-\displaystyle\sum\limits_{n=5}^\infty \frac{1}{61 \cdot f(n)}-\displaystyle\sum\limits_{n=6}^\infty \frac{1}{91 \cdot f(n)}\\+\displaystyle\sum\limits_{n=2}^\infty \frac{1}{7\cdot19 \cdot f(n)}+\displaystyle\sum\limits_{n=2}^\infty \frac{1}{7 \cdot37\cdot f(n)}+\displaystyle\sum\limits_{n=3}^\infty \frac{1}{19 \cdot37\cdot f(n)}-\displaystyle\sum\limits_{n=2}^\infty \frac{1}{7 \cdot 19\cdot 37\cdot f(3)}+\displaystyle\sum\limits_{n=2}^\infty \frac{1}{7\cdot 61 \cdot f(n)}+\displaystyle\sum\limits_{n=3}^\infty \frac{1}{19 \cdot61\cdot f(n)}\\+\displaystyle\sum\limits_{n=4}^\infty \frac{1}{37 \cdot 61\cdot f(n)}-\displaystyle\sum\limits_{n=2}^\infty \frac{1}{7\cdot19 \cdot61\cdot f(n)}-\displaystyle\sum\limits_{n=2}^\infty \frac{1}{7 \cdot37\cdot61\cdot f(n)}-\displaystyle\sum\limits_{n=3}^\infty \frac{1}{19 \cdot37\cdot61\cdot f(n)}+\displaystyle\sum\limits_{n=2}^\infty \frac{1}{7 \cdot 19\cdot 37 \cdot 61 \cdot f(n)}\\+\displaystyle\sum\limits_{n=2}^\infty \frac{1}{7 \cdot 91 \cdot f(n)}+\displaystyle\sum\limits_{n=3}^\infty \frac{1}{19 \cdot91 \cdot f(n)}+\displaystyle\sum\limits_{n=4}^\infty \frac{1}{37 \cdot 91\cdot f(n)}+\displaystyle\sum\limits_{n=5}^\infty \frac{1}{61 \cdot91 \cdot f(n)}-\displaystyle\sum\limits_{n=2}^\infty \frac{1}{7\cdot19 \cdot 91\cdot f(n)}\\-\displaystyle\sum\limits_{n=2}^\infty \frac{1}{7 \cdot37\cdot91\cdot f(n)}-\displaystyle\sum\limits_{n=3}^\infty \frac{1}{19 \cdot37\cdot91 \cdot f(n)}+\displaystyle\sum\limits_{n=2}^\infty \frac{1}{7 \cdot 19\cdot 37 \cdot 91 \cdot f(n)}-\displaystyle\sum\limits_{n=2}^\infty \frac{1}{7 \cdot 61 \cdot 91 \cdot f(n)}-\displaystyle\sum\limits_{n=3}^\infty \frac{1}{19 \cdot 61 \cdot91 \cdot f(n)}\\-\displaystyle\sum\limits_{n=4}^\infty \frac{1}{37 \cdot 61 \cdot 91 \cdot f(n)}+\displaystyle\sum\limits_{n=2}^\infty \frac{1}{7 \cdot 19 \cdot 61\cdot 91 \cdot f(n)}+\displaystyle\sum\limits_{n=2}^\infty \frac{1}{7 \cdot 37 \cdot 61 \cdot 91 \cdot f(n)}+\displaystyle\sum\limits_{n=3}^\infty \frac{1}{19 \cdot 37 \cdot 61 \cdot 91 \cdot f(n)}\\-\displaystyle\sum\limits_{n=2}^\infty \frac{1}{7 \cdot 19 \cdot 37 \cdot 61 \cdot 91 \cdot f(n)}.$$
\end{multline*}\\

\noindent We continue in this fashion until all the fraction terms on the right-hand side of the Z function in Eq. 4 are eliminated.  This, of course, will only occur at $\infty$.  But as we continue to eliminate the terms on the right-hand side of the Zeta function in Eq. 4, the value of the M-series function

$$\opM\limits_{i=2}^\infty \Big[(-1)^{i-1}\Big]\displaystyle\sum{}^i$$

\noindent which adds to unity in the numerator on the right-hand side of Eq. 3 (in this case) increases between the bounds of $-1$ and $0$ as the terms, both composite and prime, drop out of the right-hand side of Eq. 4.  To evaluate this series for the operations we have performed thus far on the Z function for the polynomial expression $f(n)=n^2-3n+1$, we rearrange the summation terms to obtain

\begin{multline*}
$$\displaystyle\sum\limits_{n=2}^\infty \frac{1}{f(n)}\cdot \Big(\color{red}-\frac{1}{7}\color{black}+\frac{1}{7\cdot 19}+\frac{1}{7\cdot 37}-\frac{1}{7\cdot 19\cdot37}+\frac{1}{7\cdot 61}-\frac{1}{7\cdot 19\cdot61}-\frac{1}{7\cdot 37\cdot61}+\frac{1}{7\cdot 19\cdot 37\cdot 61}+\frac{1}{7\cdot91}-\frac{1}{7\cdot19\cdot91}-\\\frac{1}{7\cdot37\cdot91}+\frac{1}{7\cdot19\cdot37\cdot91}-\frac{1}{7\cdot61\cdot91}+\frac{1}{7\cdot19\cdot61\cdot91}+\frac{1}{7\cdot37\cdot61\cdot91}-\frac{1}{7\cdot19\cdot37\cdot61\cdot91}...\Big)$$
\end{multline*}

\begin{flalign*}
&\quad \displaystyle\sum\limits_{n=3}^\infty \frac{1}{f(n)} \cdot\Big(\color{red}-\frac{1}{19}\color{black}+\frac{1}{19\cdot37}+\frac{1}{19\cdot61}-\frac{1}{19\cdot37\cdot61}+\frac{1}{19\cdot91}-\frac{1}{19\cdot 37\cdot91}-\frac{1}{19\cdot61\cdot91}+\frac{1}{19\cdot37\cdot 61\cdot91}...\Big)&
\end{flalign*}

\begin{flalign*}
&\quad \displaystyle\sum\limits_{n=4}^\infty \frac{1}{f(n)} \cdot\Big(\color{red}-\frac{1}{37}\color{black}+\frac{1}{37\cdot61}+\frac{1}{37\cdot91}-\frac{1}{37\cdot61\cdot91}...\Big)&
\end{flalign*}

\begin{flalign*}
&\quad \displaystyle\sum\limits_{n=5}^\infty \frac{1}{f(n)}\cdot \Big(\color{red}-\frac{1}{61}\color{black}+\frac{1}{61\cdot91}...\Big)&
\end{flalign*}

\begin{flalign*}
&\quad \displaystyle\sum\limits_{n=6}^\infty \frac{1}{f(n)}\cdot \Big(\color{red}-\frac{1}{91}\color{black}...\Big)&
\end{flalign*}

\begin{flalign*}
&\quad \displaystyle\sum\limits_{n=7}^\infty \frac{1}{f(n)}.&
\end{flalign*}

\noindent Once the summation terms are thus grouped in this fashion according to their lower limits, we begin collecting like terms across the groups to organize the sums according to the number of terms which their product will contain when they are expanded with their coefficients.  Referring to the grouping above, we begin by extracting the summations with the fewest terms which have been highlighted in \color{red} red \color{black}, e.g.,

\begin{flalign*}
&\quad \color{red}-\frac{1}{7}\color{black}\cdot\displaystyle\sum\limits_{n=2}^\infty \frac{1}{f(n)}\color{red}-\frac{1}{19}\color{black}\cdot\displaystyle\sum\limits_{n=3}^\infty \frac{1}{f(n)}\color{red}-\frac{1}{37}\color{black}\cdot\displaystyle\sum\limits_{n=4}^\infty \frac{1}{f(n)}\color{red}-\frac{1}{61}\color{black}\cdot\displaystyle\sum\limits_{n=5}^\infty \frac{1}{f(n)}\color{red}-\frac{1}{91}\color{black}\cdot\displaystyle\sum\limits_{n=6}^\infty \frac{1}{f(n)}\color{red}-\color{black}\cdots=&
\end{flalign*}

\begin{flalign*}
&\quad \color{red}-\frac{1}{7}\color{black}\cdot\Big(\frac{1}{7}+\frac{1}{19}+\frac{1}{37}+\frac{1}{61}+\frac{1}{91}+\cdots\Big)&
\end{flalign*}

\begin{flalign*}
&\quad \color{red}-\frac{1}{19}\color{black}\cdot\Big(\frac{1}{19}+\frac{1}{37}+\frac{1}{61}+\frac{1}{91}+\cdots\Big)&
\end{flalign*}

\begin{flalign*}
&\quad \color{red}-\frac{1}{37}\color{black}\cdot\Big(\frac{1}{37}+\frac{1}{61}+\frac{1}{91}+\cdots\Big)&
\end{flalign*}

\begin{flalign*}
&\quad \color{red}-\frac{1}{61}\color{black}\cdot\Big(\frac{1}{61}+\frac{1}{91}+\cdots\Big)&
\end{flalign*}

\begin{flalign*}
&\quad \color{red}-\frac{1}{91}\color{black}\cdot\Big(\frac{1}{91}+\cdots\Big)&
\end{flalign*}

\begin{flalign*}
&\quad \color{red}-\color{black}\cdots&
\end{flalign*}

\noindent which yields the first term in the infinite series, or

$$\opM\limits_{i=2}^{} \Big[(-1)^{i-1}\Big]\displaystyle\sum{}^i$$

\noindent previously defined as

$$-\Sigma_1 \Sigma_2 = -\sum\limits_{i=2}^\infty \sum\limits_{j=i}^\infty \frac{1}{f(i)} \frac{1}{f(j)}.$$

\noindent The second term of the M-series of Eq. 3 also formed by collecting like terms (highlighted in \color{red} red) \color{black} across the original grouping as follows:

\begin{multline*}
$$\displaystyle\sum\limits_{n=2}^\infty \frac{1}{f(n)}\cdot \Big(-\frac{1}{7}\color{red}+\frac{1}{7\cdot 19}+\frac{1}{7\cdot 37}\color{black}-\frac{1}{7\cdot 19\cdot37}\color{red}+\frac{1}{7\cdot 61}\color{black}-\frac{1}{7\cdot 19\cdot61}-\frac{1}{7\cdot 37\cdot61}+\frac{1}{7\cdot 19\cdot 37\cdot 61}\color{red}+\frac{1}{7\cdot91}\color{black}-\frac{1}{7\cdot19\cdot91}-\\\frac{1}{7\cdot37\cdot91}+\frac{1}{7\cdot19\cdot37\cdot91}-\frac{1}{7\cdot61\cdot91}+\frac{1}{7\cdot19\cdot61\cdot91}+\frac{1}{7\cdot37\cdot61\cdot91}-\frac{1}{7\cdot19\cdot37\cdot61\cdot91}...\Big)$$
\end{multline*}

\begin{flalign*}
&\quad \displaystyle\sum\limits_{n=3}^\infty \frac{1}{f(n)} \cdot\Big(-\frac{1}{19}\color{red}+\frac{1}{19\cdot37}+\frac{1}{19\cdot61}\color{black}-\frac{1}{19\cdot37\cdot61}\color{red}+\frac{1}{19\cdot91}\color{black}-\frac{1}{19\cdot 37\cdot91}-\frac{1}{19\cdot61\cdot91}+\frac{1}{19\cdot37\cdot 61\cdot91}...\Big)&
\end{flalign*}

\begin{flalign*}
&\quad \displaystyle\sum\limits_{n=4}^\infty \frac{1}{f(n)}\cdot\Big(-\frac{1}{37}\color{red}+\frac{1}{37\cdot61}+\frac{1}{37\cdot91}\color{black}-\frac{1}{37\cdot61\cdot91}...\Big)&
\end{flalign*}

\begin{flalign*}
&\quad \displaystyle\sum\limits_{n=5}^\infty\frac{1}{f(n)}\cdot \Big(-\frac{1}{61}\color{red}+\frac{1}{61\cdot91}\color{black}...\Big)&
\end{flalign*}

\begin{flalign*}
&\quad \displaystyle\sum\limits_{n=6}^\infty \frac{1}{f(n)}\cdot \Big(-\frac{1}{91}...\Big)&
\end{flalign*}

\begin{flalign*}
&\quad \displaystyle\sum\limits_{n=7}^\infty \frac{1}{f(n)}&
\end{flalign*}

\noindent And these terms can be extracted from the original grouping and organized in the same fashion:

\begin{flalign*}
&\quad\Big(\color{red}\frac{1}{7\cdot19}+\frac{1}{7\cdot37}+\frac{1}{7\cdot61}+\frac{1}{7\cdot91}+\cdots\color{black}\Big)\cdot\displaystyle\sum\limits_{n=2}^\infty \frac{1}{f(n)}+&
\end{flalign*}

\begin{flalign*}
&\quad\Big(\color{red}\frac{1}{19\cdot37}+\frac{1}{19\cdot61}+\frac{1}{19\cdot91}+\cdots\color{black}\Big)\cdot\displaystyle\sum\limits_{n=3}^\infty \frac{1}{f(n)}+&
\end{flalign*}

\begin{flalign*}
&\quad\Big(\color{red}\frac{1}{37\cdot61}+\frac{1}{37\cdot91}+\cdots\color{black}\Big)\cdot\displaystyle\sum\limits_{n=4}^\infty \frac{1}{f(n)}+&
\end{flalign*}

\begin{flalign*}
&\quad\Big(\color{red}\frac{1}{61\cdot91}+\cdots\color{black}\Big)\cdot\displaystyle\sum\limits_{n=5}^\infty \frac{1}{f(n)}+\cdots&
\end{flalign*}

\begin{flalign*}
&\quad=&
\end{flalign*}

\begin{flalign*}
&\quad \Big(\color{red}\frac{1}{7\cdot19}+\frac{1}{7\cdot37}+\frac{1}{7\cdot61}+\frac{1}{7\cdot91}+\cdots\color{black}\Big)\cdot\Big(\frac{1}{7}+\frac{1}{19}+\frac{1}{37}+\frac{1}{61}+\frac{1}{91}+\cdots\Big)+&
\end{flalign*}

\begin{flalign*}
&\quad \Big(\color{red}\frac{1}{19\cdot37}+\frac{1}{19\cdot61}+\frac{1}{19\cdot91}+\cdots\color{black}\Big)\cdot\Big(\frac{1}{19}+\frac{1}{37}+\frac{1}{61}+\frac{1}{91}+\cdots\Big)+&
\end{flalign*}

\begin{flalign*}
&\quad \Big(\color{red}\frac{1}{37\cdot61}+\frac{1}{37\cdot91}+\cdots\color{black}\Big)\cdot\Big(\frac{1}{37}+\frac{1}{61}+\frac{1}{91}+\cdots\Big)+&
\end{flalign*}

\begin{flalign*}
&\quad \Big(\color{red}\frac{1}{61\cdot91}+\cdots \color{black}\Big)\cdot\Big(\frac{1}{61}+\frac{1}{91}+\cdots\Big)+\cdots&
\end{flalign*}

\noindent which yields the second term in the infinite series, or

$$\opM\limits_{i=3}^{} \Big[(-1)^{i-1}\Big]\displaystyle\sum{}^i$$

\noindent previously defined as

$$+\Sigma_1 \Sigma_2 \Sigma_3 = +\sum\limits_{i=2}^\infty \sum\limits_{j=i}^\infty \sum\limits_{k=j+1}^\infty \frac{1}{f(i)} \frac{1}{f(j)} \frac{1}{f(k)}.$$

\noindent The third term of the M-series of Eq. 3 is formed by collecting the next level of terms across the summations in the original grouping:\\

\begin{multline*}
$$\displaystyle\sum\limits_{n=2}^\infty \frac{1}{f(n)}\cdot \Big(-\frac{1}{7}+\frac{1}{7\cdot 19}+\frac{1}{7\cdot 37}\color{red}-\frac{1}{7\cdot 19\cdot37}\color{black}+\frac{1}{7\cdot 61}\color{red}-\frac{1}{7\cdot 19\cdot61}-\frac{1}{7\cdot 37\cdot61}\color{black}+\frac{1}{7\cdot 19\cdot 37\cdot 61}+\frac{1}{7\cdot91}\color{red}-\frac{1}{7\cdot19\cdot91}-\color{black}\\\color{red}\frac{1}{7\cdot37\cdot91}\color{black}+\frac{1}{7\cdot19\cdot37\cdot91}\color{red}-\frac{1}{7\cdot61\cdot91}\color{black}+\frac{1}{7\cdot19\cdot61\cdot91}+\frac{1}{7\cdot37\cdot61\cdot91}-\frac{1}{7\cdot19\cdot37\cdot61\cdot91}...\Big)$$
\end{multline*}

\begin{flalign*}
&\quad \displaystyle\sum\limits_{n=3}^\infty \frac{1}{f(n)} \cdot\Big(-\frac{1}{19}+\frac{1}{19\cdot37}+\frac{1}{19\cdot61}\color{red}-\frac{1}{19\cdot37\cdot61}\color{black}+\frac{1}{19\cdot91}\color{red}-\frac{1}{19\cdot37\cdot91}-\frac{1}{19\cdot61\cdot91}\color{black}+\frac{1}{19\cdot37\cdot 61\cdot91}...\Big)&
\end{flalign*}

\begin{flalign*}
&\quad \displaystyle\sum\limits_{n=4}^\infty \frac{1}{f(n)}\cdot\Big(-\frac{1}{37}+\frac{1}{37\cdot61}+\frac{1}{37\cdot91}\color{red}-\frac{1}{37\cdot61\cdot91}\color{black}...\Big)&
\end{flalign*}

\begin{flalign*}
&\quad \displaystyle\sum\limits_{n=5}^\infty\frac{1}{f(n)}\cdot \Big(-\frac{1}{61}+\frac{1}{61\cdot91}...\Big)&
\end{flalign*}

\begin{flalign*}
&\quad \displaystyle\sum\limits_{n=6}^\infty \frac{1}{f(n)}\cdot \Big(-\frac{1}{91}...\Big)&
\end{flalign*}

\begin{flalign*}
&\quad \displaystyle\sum\limits_{n=7}^\infty \frac{1}{f(n)}&
\end{flalign*}

\noindent and are extracted in the same fashion:

\begin{flalign*}
&\quad\Big(\color{red}-\frac{1}{7\cdot19\cdot37}-\frac{1}{7\cdot19\cdot61}-\frac{1}{7\cdot37\cdot61}-\frac{1}{7\cdot19\cdot91}-\frac{1}{7\cdot37\cdot91}-\frac{1}{7\cdot61\cdot91}-\cdots\color{black}\Big)\cdot\displaystyle\sum\limits_{n=2}^\infty \frac{1}{f(n)}+&
\end{flalign*}

\begin{flalign*}
&\quad\Big(\color{red}-\frac{1}{19\cdot37\cdot61}-\frac{1}{19\cdot37\cdot91}-\frac{1}{19\cdot61\cdot91}-\cdots\color{black}\Big)\cdot\displaystyle\sum\limits_{n=3}^\infty \frac{1}{f(n)}+&
\end{flalign*}

\begin{flalign*}
&\quad\Big(\color{red}-\frac{1}{37\cdot61\cdot91}-\cdots\color{black}\Big)\cdot\displaystyle\sum\limits_{n=4}^\infty \frac{1}{f(n)}+\cdots&
\end{flalign*}

\begin{flalign*}
&\quad=&
\end{flalign*}

\begin{flalign*}
&\quad \Big(\color{red}-\frac{1}{7\cdot19\cdot37}-\frac{1}{7\cdot19\cdot61}-\frac{1}{7\cdot37\cdot61}-\frac{1}{7\cdot19\cdot91}-\frac{1}{7\cdot37\cdot91}-\frac{1}{7\cdot61\cdot91}-\cdots\color{black}\Big)\cdot\Big(\frac{1}{7}+\frac{1}{19}+\frac{1}{37}+\frac{1}{61}+\frac{1}{91}+\cdots\Big)+&
\end{flalign*}

\begin{flalign*}
&\quad \Big(\color{red}-\frac{1}{19\cdot37\cdot61}-\frac{1}{19\cdot37\cdot91}-\frac{1}{19\cdot61\cdot91}-\cdots\color{black}\Big)\cdot\Big(\frac{1}{19}+\frac{1}{37}+\frac{1}{61}+\frac{1}{91}+\cdots\Big)+&
\end{flalign*}

\begin{flalign*}
&\quad \Big(\color{red}-\frac{1}{37\cdot61\cdot91}\color{black}\Big)\cdot\Big(\frac{1}{37}+\frac{1}{61}+\frac{1}{91}+\cdots\Big)+\cdots&
\end{flalign*}\\

\noindent which yields the third term in the infinite series, or

$$\opM\limits_{i=4}^{} \Big[(-1)^{i-1}\Big]\displaystyle\sum{}^i$$

\noindent previously defined as

$$-\Sigma_1 \Sigma_2 \Sigma_3 \Sigma_4 = -\sum\limits_{i=2}^\infty \sum\limits_{j=i}^\infty \sum\limits_{k=j+1}^\infty \sum\limits_{l=k+1}^\infty \frac{1}{f(i)} \frac{1}{f(j)} \frac{1}{f(k)} \frac{1}{f(l)}.$$

\noindent We now see that a pattern emerges in which it is clear that an infinite sum of infinite sums will be obtained if one continues to eliminate terms from the big Zeta function in Eq. 4 using the non-sieving approach.  It is hypothesized that once the infinity of terms are sieved out of the right hand side of the Z function in Eq. 4, a constant will emerge for the M-series in Eq. 3 which takes on a finite value between

$$\color{red}-1\le\color{black} \quad \opM\limits_{i=2}^\infty \Big[(-1)^{i-1}\Big]\displaystyle\sum{}^i \quad \color{red}\le{0}\color{black}$$

\noindent and which may tell us something about the density of prime numbers for an integer-valued polynomial function $f(n)$ for some limit of the domain $n \le x$.\\

\noindent As implied in row $i$ of Fig. 1, if one applies the realm of positive integers $i \ge 2$ to the Euler zeta function using the non-sieving approach so that the product term in the denominator of the right-hand side of Eq. 1 includes all the prime \underline{and} composite terms, then the numerator on the right-hand side of Eq. 1 will tend to zero as the limit $x$ tends to infinity.  Thus in this case, one could associate the cumulative value of all the composite terms in the product of the denominator of the right-hand side of the Euler zeta function in Eq. 1 to an M-series value of $-1$ in the numerator of the right-hand side of Eq. 3, because when all the composite terms are included in the product term in the denominator of the right-hand side of the Euler zeta function in Eq. 1, the numerator on the right-hand side of the big Zeta function in Eq. 3 will tend to zero as $x$ approaches $\infty$ (i.e., $1+(M)=0$).  Compare that to the value of unity obtained in the numerator of the right-hand side of Eq. 1 when all the composite terms are sieved out of the product term in the denominator of the right-hand side of the equation as Euler did when he discovered this famous relationship between the composite and prime numbers.  So one might postulate that the total effect of sieving the composite terms out of the product term in the denominator of the right-hand side of the Euler zeta function in Eq. 1 is that the numerator in the right-hand side of the equation increases from zero to unity as the composite terms are eliminated from the product term as $x$ approaches $\infty$.  It was observed that while the M-series

$$\displaystyle\opM\limits_{i=2}^\infty \Big{[(-1)^{i-1}\Big]\displaystyle\sum{}^i}$$

\noindent represents an infinite sum of infinite sums, its value is constrained between $-1$ and $0$ throughout the application of the non-sieving function to the realm of positive integers as well as to the range of integer calculations yielded by the integer-valued polynomial function demonstrated in the aforementioned example.\\

\begin{center}
\includegraphics[scale=0.6]{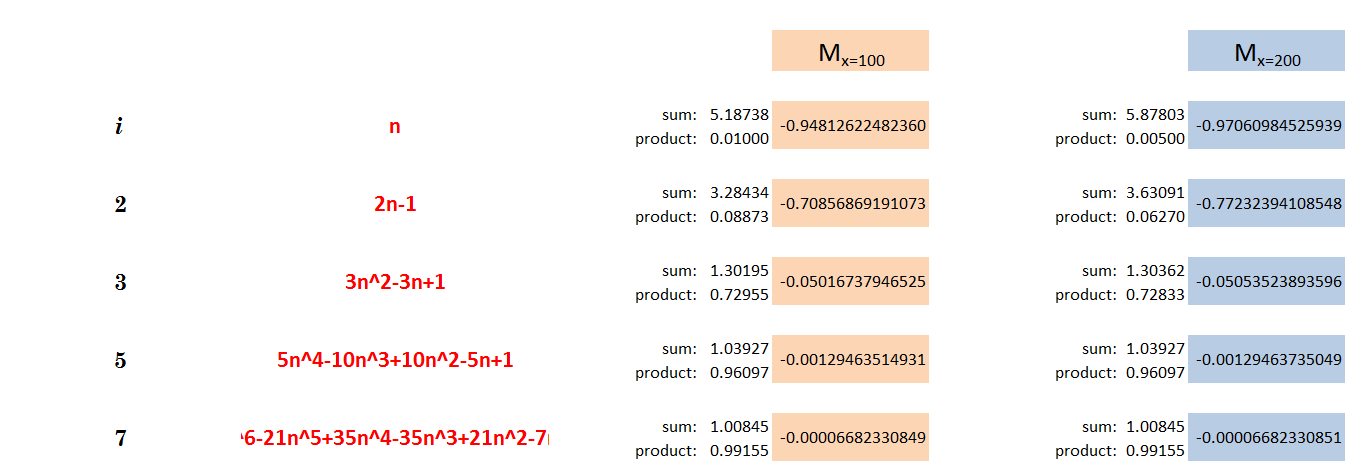}
\end{center}

\begin{center}
\textit{Fig. 1 - Values of the M-series for the realm of the positive integers $i$ and for the first four prime powers $p=2,3,5,7$ of the prime shell function for the limits $x=100$ and $x=200$}
\end{center}

\noindent Descending rows below row $i$ in Fig. 1 tabulate the calculations of $M_{x=100}$ and $M_{x=200}$, respectively, that result from the non-sieving application of the zeta function to the prime shell function $$n^p-(n-1)^p$$ with prime powers $2, 3, 5$ and $7$, respectively.  Included next to each calculation of $M_{x=100}$ and $M_{x=200}$ are the $Ms$ (sum) and $\Pi$ (product) calculations for that row of the prime shell function which were used to calculate the values of $M_{x=100}$ and $M_{x=200}$.  Fig. 2 is a graph plotted to superimpose the trends of $M_{x=100}$ and $M_{x=200}$ in Fig. 1 for the case of the non-sieving application of the Euler zeta function and illustrates the change in the value of $M_{x=100}$ and $M_{x=200}$ as the prime power in the prime shell function increases from $2$ to $7$ and as the range of the base $n$ of the prime shell function increases from $x=100$ to $x=200$.  Notice the trend is that the M-series value increases toward $0$ as the prime power in the prime shell equation increases but that the change in the M-series values due to the increase of the range $x$ from $100$ to $200$ is very slight.  It is anticipated that as $x$ approaches $\infty$ for row $i$ in the case of the application of the non-sieving zeta function to the realm of positive integers that the value of $M$ will approach $-1$.  It is also anticipated that the limit of $M$ in the subsequent rows of the application of the non-sieving approach to the prime shell functions with powers $2, 3, 5$ and $7$ will approach a constant value as the range of $x$ of the base $n$ approaches $\infty$.  From the trends manifested by the M-series values in Figures 1 and 2, these values are estimated to be close to

$$M_{{i}_{(x=\infty)}}\approx -1$$ $$M_{{2}_{(x=\infty)}}\approx -0.8$$ $$M_{{3}_{(x=\infty)}}\approx -0.051$$ $$M_{{5}_{(x=\infty)}}\approx -0.0013$$ $$M_{{7}_{(x=\infty)}}\approx -0.000067.$$

\begin{center}
\includegraphics[scale=1.1]{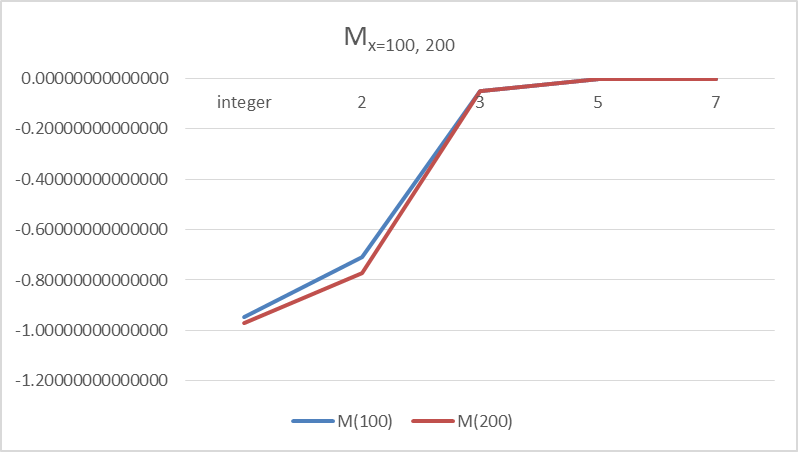}
\end{center}
\begin{center}

\textit{Fig. 2 - Graph of $M$ superimposed for $x=100$ and $x=200$ for \textit{i}$=$integer and for the first four prime powers $p=2,3,5,7$ of the prime shell polynomial functions.  The M-series is asymptotic to $0$ as the prime power $p$ approaches $\infty$} \\
\end{center}

\noindent Table 1 reveals the M-series values of $M_{x=100}$ and $M_{x=200}$ for the non-sieving application of the Euler zeta function to the integer domain $2\le i \le 100$ and $2\le i \le 200$.  It is seen that as the limit $x$ increases from $x=100$ to $x=200$ that the value of the M-series approaches $-1$.  Table 2 includes the calculated values of $M_{x=100}$ and $M_{x=200}$ for the prime shell function $n^p-(n-1)^p$ for prime powers $2, 3, 5$ and $7$.  It is anticipated that these M-series values, as they are refined by increasing the limit $x$ of the domain of the base $n$ of the function, will be able to tell us something about the density of prime numbers yielded by integer-valued polynomial functions within the range of the domain as base $n$ is increased to some limit $x$ for the prime shell functions.  Further study should reveal more about these M-series properties. \\

\begin{center}

\begin{tabular}{|c||c|c|c||c|c|c|}
\hline
prime number $\ge2$ & $\pi \le 100$ & $Ls \le 100$ & ${\displaystyle\opM\limits_{}^{}}_{x=100}$ &  $\pi \le 200$ & $Ls \le 200$ & ${\displaystyle\opM\limits_{}^{}}_{x=200}$ \\
\hline
\textit{i} & 25 & 29.99144 & -0.94812622482360 & 46 & 50.04329 & -0.97060984525939 \\
\hline
\end{tabular}
\end{center}

\begin{center}
\textit{Table 1}\\
\end{center}

\begin{center}

\footnotesize

\begin{tabular}{|c||c|c|c||c|c|c|}

\hline
prime power $p$ & $\Pi \le f(100)$ & $Ms \le f(100)$ & ${\displaystyle\opM\limits_{}^{}}_{x=100}$ &  $\Pi \le f(200)$ & $Ms \le f(200)$ & ${\displaystyle\opM\limits_{}^{}}_{x=200}$ \\
\hline
2 & 44 & 42.75969 & -0.70856869191073 & 76 & 78.48273 & -0.77232394108548 \\
\hline
3 & 43 & 29.01307 & -0.05016737946525 & 72 & 53.06455 & -0.05053523893596 \\
\hline
5 & 18 & 19.71488 & -0.00129463514931 & 32 & 35.92022 &-0.00129463735049 \\
\hline
7 & 24 & 15.71077 & -0.00006682330849 & 40 & 28.56513 & -0.00006682330851 \\
\hline
\end{tabular}
\end{center}

\begin{center}
\textit{Table 2}\\
\end{center}

\noindent It is further anticipated by the author that the finite values yielded by the M-series in non-sieving applications of the Euler zeta function in the case of integer-valued polynomials in general, within the ranges that those polynomials are evaluated, may tell us something about the density of numbers within the ranges of those functions compared to the number of primes there are on the real number line less than or equal to some upper limit $x$ of the domain of base $n$ that served as inputs to those integer-valued polynomials.  More time will be needed to evaluate any such possible relationships that may exist.  Until such time, the author proposes the following theorem:

\begin{theorem}
\noindent When a non-sieving application of the Euler zeta function is applied to process values generated by an integer-valued polynomial, then there is an infinite series

$$\displaystyle\opM\limits_{i=2}^\infty \Big{[(-1)^{i-1}\Big]\displaystyle\sum{}^i}$$

\noindent that arises which adds to unity in the numerator of the product term in the zeta function to make both sides of the equation equal.  The value of this M-series is bounded by $-1$ and $0$ in the non-sieving application of the Euler zeta function.
\end{theorem}

\noindent It is the author's hope that the results of this study will open the doors for further research into prime number frequency among the number sequences yielded by integer-valued polynomials which use as their domain the realm of positive integers, as well as other integer-valued polynomial functions in general which may use other domains which enable those functions to generate prime numbers in a seemingly random fashion among the composite numbers within its range of calculations.  It is also the author's hope that the results of this research will present some open-ended problems for which mathematicians and mathematics students will want to ponder and tackle in attempt to shed more light on prime number theory in general.

\end{document}